\input amstex
\documentstyle{amsppt}
\pageheight{50.5pc}
\pagewidth{32pc}
\define\1{\hbox{\rm 1}\!\! @,@,@,@,@,\hbox{\rm I}}

\topmatter
\title
{
Limit theorems for rarefaction of set of diffusion processes by boundaries.
\footnotemark{}
}
\endtitle
\footnotetext{This research was supported (in part) by the
Ministry of Education and Science of Ukraine, project
No 01.07/103 and University Salerno, Italy.}
\author
{Aniello Fedullo and Vitalii A. Gasanenko}
\endauthor

\address
{Universita Degli Studi Di Salerno, Via Ponte Don Melillo, 84084
 Fisciano (SA) Italia}
\endaddress
\email
{afedullo\@unisa.it}
\endemail
\address
{Institute of Mathematics,National Academy of
Science of Ukraine, Tereshchenkivska 3, 01601, Kiev, Ukraine}
\endaddress
\email
{gs\@imath.kiev.ua}
\endemail
\subjclass
{60 J 60}
\endsubjclass
\keywords
{Stochastic differential equations,  solution of parabolic equations,
 eigenvalues problem, generating function}
\endkeywords
\subjclass
{60 J 60}
\endsubjclass
\email
{afedulo\@unisa.it}
\endemail
\email
{gs\@imath.kiev.ua}
\endemail
\abstract
This paper is devoted to the study of the following problem.
We have set of diffusion  processes  with absorption on boundaries
in some region  at initial time $t=0$. It is required to estimate of number
of the  unabsorbed processes for the fixed time ~$\tau>0$. The number of initial processes
is considered as function of $\tau$ and $\tau\to\infty$.
\endabstract

\endtopmatter

\document

Consider the ~$N$~ independent diffusion processes
which start from different points ~$x_{k}\in Q,\quad k=\overline{1,N}$~~
in the initial moment of time ~$t=0$.
The region  ~~$ Q\subset R^{d}$~~ is open connected region and it is limited by
smooth surface $\partial Q$. All processes are diffusion processes
with absorption on the boundary ~$\partial Q$.
 These processes are solutions
of the following stochastic differential  equations in $Q$
with absorption

$$
d\xi(t)=a(t,\xi(t))dt + \sum\limits_{i=1}^{d}b_{i}(t,\xi (t))dw^{(k)}_{i}(t)
\quad \xi(t)\in R^{d}\eqno(1)
$$

$$
 b_{i}(t,x),~ a(t,x): R_{+}\times R^{d}\to R^{d}.
$$

with an initial condition: ~$\xi(0)=\gamma_{k}\in Q.$

Here the ~$W^{(k)}(t)=(w_{i}^{(k)}(t),\quad 1\leq i\leq d),~~1\leq k\leq N$~
 are independent in totality $d$- dimensional Wiener processes.

Thus, these processes have the identical diffusion matrices and
shift vectors , but they have different initial states.

We will be interested by distribution of number of processes yet
not absorbed by boundary ~$\partial Q$~ in the moment of time ~$\tau$.

This task was offered in [1] as the mathematical model of cleaning of gas
from particles (dust, microbes and ect.).
Cleaning consists in pass of gas
 with speed ~$v$~ through the pipe of length ~$l$~ and
with permanent section ~$Q\cup \partial Q$~.
Walls of pipe absorb the particles.
The estimation of efficiency of cleaning can be reduced
to the solution of this problem at ~$\tau=l/v$. In supposition,
that  particles evolve independently from each other and
their movement can be presented as the solutions of (1).

We make the following assumptions:
a number and positions of particles
 are defined  by a given determinate limited measure
~$N(B,\tau)$ at the initial time.
Thus ~$N(B,\tau)$~ is equal to number of points ~$x_{k}$~ in the set ~$B$~
and ~$N=N(Q,\tau)<\infty$ for fixed $\tau>0$.

Note, the case when ~$N(\cdot,\tau)=N(\cdot)$~ there is random Poisson measure on a circle

~$C=\{(x,y):\quad x^{2}+y^{2}\leq r^{2}\}$ was considered in [1].

There is exact formula of distribution function for number of remaining processes
 in this article. However, this formula consists diffucult computed
functions.
The present article is devoted to case when initial number of processes
$N(Q,\tau)$ depends on the final time $\tau$:
$N(Q,\tau)\to\infty$ when $\tau\to\infty$.
We shall obtain conditions of such dependence which leads
to simple limit distribution function of number not absorbed processes.

Suppose that the region $Q$ is bounded and boundary $\partial Q$ is  Lyapunov surface
~$C^{(1,\lambda)}$ [2].
We will consider the following case

$a(t,x)=(\underbrace{0,\dots,0}_{d}),\quad b_{i}(t,x)=b_{i}=
(b_{i1},\dots,b_{id}),~~1\leq i\leq d.$

We will define  matrix
 ~$\sigma=B^{T}B,\quad B=(b_{ij}),~1\leq i,j\leq d\quad
\sigma=(\sigma_{ij}),1\leq i,j\leq d$ and differential operator
$A:\sum\limits_{1\leq i,j\leq d}\sigma_{ij}\frac{\partial^{2}}
{\partial x_{i}\partial x_{j}}.$
Let $\sigma$ be a matrix with  the following property

$$
\sum\limits_{1\leq i,j\leq d}\sigma_{ij}z_{i}z_{j}
\geq \mu |\vec z|^{2}.
$$

Here ~$\mu$,~ there is  fixed positive number, and
$\vec z=(z_{1},\cdots, z_{d})$~ there is an arbitrary real vector.

This operator acts in the following space

$$
H_{A}=\{u: u\in L_{2}(Q)\cap Au\in L_{2}(Q)\cap u(\partial Q)=0\}
$$

with inner product $(u,v)_{A}=(Au,v)$.
Here $(,)$ is inner product in $L_{2}(Q)$.
 The operator ~$A$~ is  positive
operator [2]. It is known that the following eigenvalues problem

$$
Au=-\lambda u,\quad u(\partial Q)=0
$$

has infinity set of real eigenvalues  $\lambda_{i}\to\infty$ and

$$
0<\lambda_{1}<\lambda_{2}<\cdots<\lambda_{s}<\cdots.
$$

The corresponding eigenfunctions

$$
f_{11},\dots,f_{1n_{1}},\cdots,f_{s1},\dots,f_{sn_{s}},\cdots
$$
form complete system of functions both in $H_{A}$ and $L_{2}^{0}(Q):=
\{u: u\in L_{2}(Q)\cap u(\partial Q)=0\}$. Here number $n_{k}$ is equal to
multiplicity  of eigenvalue  $\lambda_{k}$.

We will denote by ~$\eta(\tau)$~ the number of remaining processes in the region ~$D$~
at the moment ~$\tau$.

We will also assume that $\sigma$-additive  measure $\nu$ is given on the
$\Sigma_{\nu}$- algebra sets from $Q,\quad \nu(Q)<\infty.$
All eigenfunctions
$f_{ij}:Q\to R^{1}$ are $(\Sigma_{\nu},\Sigma_{Y})$ measurable. Here
$\Sigma_{Y}$ is system of Borel sets from $R^{1}$.
Let $\Rightarrow$ denote the  weak convergence of
random values or measures.

Let us denote by ~$\nu_{\tau}(\cdot)$ the measure

$$
\nu_{\tau}(B)=\exp(-\frac{\tau}{2}\lambda_{1})N(B,\tau).
$$

where $B\in\Sigma_{\nu}$ .

By definition of measure ~$\nu_{\tau}(\cdot)$, we have
$$
d\nu_{\tau}(x)=
\cases
\exp(-\frac{\tau}{2}\lambda_{1}),&  \hbox{if}\quad x=x_{k},\quad k=1,\cdots, N(Q,\tau) \\
0 ,& \hbox{otherwise}.
\endcases
$$

\proclaim
{\bf Theorem 1} Suppose $N(\cdot, \tau)$ satisfies the condition

$$
\nu_{\tau}(\cdot)\mathop{\Rightarrow}\limits_{\tau\to\infty}\nu(\cdot).
$$

Then $\eta(\tau)\Rightarrow\eta$ if $\tau\to\infty$ where $\eta$ has Poisson
distribution function with parameter $a=\int\limits_{Q} F(x)d\nu(x)$ and
$F(x)=\sum\limits_{i=1}^{n_{1}}f_{1i}(x)c_{1i},\quad
c_{1i}=\int\limits_{Q}f_{1i}(x)d x.$
\endproclaim

\demo
{\bf Proof}   Consider the following initial-boundary problem

$$
\frac{\partial u}{\partial t}=\frac{1}{2}\sum\limits_{1\leq i,j\leq d}
\sigma_{ij}\frac{\partial^{2}u}{\partial x_{i}\partial x_{j}}
,\quad x\in Q;
$$

$$
u(0,x)= 1, \quad x\in Q
$$

$$
u(t,x)=0 \quad \hbox{if}\quad x\in \partial Q, \quad t\geq 0\eqno(2)
$$

It is known [3, sec. VIII], that ~ $u(\tau,x)$~ is equal to probability to remain in the
 region ~$Q$~  at the moment ~$\tau$~ for a diffusion process which
 was in a point $x$~at the  initial moment (~$\xi(0)= x\in Q$).

Introduce indicators
$$
\chi_{k}(\tau)=
\cases
1,&  \hbox{if}\quad k-\hbox{th  particle belongs to}\quad Q
\hbox{ at the moment}\quad \tau \\
0 ,& \hbox{if}\quad k-\hbox{th particle absorbed by the moment}\quad \tau.
\endcases
$$

These indicators are mutually independence by assumption. Thus the following
relations are correct

$$
\eta(\tau)=\sum\limits_{k=1}^{N(Q,\tau)}\chi_{k}(\tau),
$$
$$
E s^{\eta(\tau)}=\prod_{k=1}^{N(Q,\tau)}E s^{\chi_{k}(\tau)}=
\prod_{k=1}^{N(Q,\tau)}(su(\tau,\gamma_{k})+(1-u(\tau,\gamma_{k}));
$$
$$
\ln E s^{\eta(\tau)}=
\sum\limits_{k=1}^{N(Q,\tau)}\ln (1-u(\tau,\gamma_{k})
(1-s)).
$$

Here $0\leq s\leq 1$. As $0\leq u(\tau,\gamma_{k})\leq 1$ then we have the following
inequality from the last

$$
|\ln E s^{\eta(\tau)} -
\sum\limits_{k=1}^{N(C,\tau)}u(\tau,\gamma_{k})(1-s)|\leq
\alpha\sum\limits_{k=1}^{N(C,\tau)}u^{2}(\tau,\gamma_{k}),\eqno(3)
$$

Here $\alpha<\infty$. Define the value of $u(\tau,\gamma_{k})$.

We shall assume
that system of functions $\{f_{ij}(x), i\geq 1, 1\leq j\leq n_{i}\}$
is orthonormalized with respect to space $L_{2}^{0}(Q)$.
The ordinary argumentaion (see, for example, [2, sec. 22])
leads to the definition of solution of problem (2) of the  form

$$
u(t,x)=\sum\limits_{k=1}^{\infty}\exp(-\frac{t}{2}\lambda_{k})
\sum\limits_{j=1}^{n_{k}}c_{kj}f_{kj}( x)
$$

where coefficients $c_{ij}$ are equal to coefficients of decomposition
of initial value (unit)  by system of functions $f_{ij}$:
$c_{ij}=\int\limits_{Q}f_{ij}(x)d x$. The Parseval - Steklov equality
is true for these coefficients:

$$
\sum\limits_{k=1}^{\infty}\sum\limits_{j=1}^{n_{k}}c_{kj}^{2}=|Q|.\eqno(4)
$$

Further

$$
\sum\limits_{k=1}^{N(Q,\tau)}u(\tau,\gamma_{k})=\int\limits_{Q}u(\tau,x)
d\nu_{\tau}(x)=
$$

$$
= \int\limits_{Q}F(x)d\nu_{\tau}(x)+\int\limits_{Q}s_{\tau}(x)d\nu_{\tau}(x),
\eqno(5)
$$

here

$$
s_{\tau}(x):=
\sum\limits_{k\geq 2}\exp\left(-\frac{\tau}{2}(\lambda_{k}-\lambda_{1})\right)
\sum\limits_{j=1}^{n_{k}}c_{kj}f_{kj}(x).
$$

As function $F(x)$ is continuous and bounded function on the $\bar Q$ then
under the condition of the theorem, we obtain [4]

$$
\int\limits_{Q}F(x)d\nu_{\tau}(x)\mathop{\longrightarrow}\limits_{\tau\to
\infty}\int\limits_{Q}F(x)d\nu(x).
$$

In order to estimate of $s_{\tau}(x)$, we will give the result
from monorgaphy [5,Thm. 17.5.3]

Consider  the sums of  eigenfunctions of the form

$$
e(x,\lambda)=\sum\limits_{\lambda_{k}\leq \lambda}\sum\limits_{j=1}^{n_{k}}
f^{2}_{kj}(x)
$$

then
$$
\sup\limits_{x\in Q}\sqrt{e(x,\lambda)}\leq C \lambda^{\frac{d}{2}}.
$$

The asymptotic  of eigenvalues ~$\lambda_{k}$~ under ~$k\to\infty$~
is defined by the following eniqualities [2, sec. 18]

$$
c_{1}k^\frac{2}{d}\leq \lambda_{k}\leq c_{2} k^\frac{2}{d}, \quad
\hbox{where}\quad c_{1},~c_{2}=const.
$$

The last and (4) and  Caushy-Bunyakovskii inequality
lead to the following convergence under
$\tau\to\infty$

$$
|s_{\tau}( x)|\leq
\sum\limits_{k\geq 2}\exp\left(-\frac{\tau}{2}(\lambda_{k}-\lambda_{1})\right)
\sqrt{\sum\limits_{j=1}^{n_{k}}
c^{2}_{kj}}\sqrt{\sum\limits_{m=1}^{n_{k}}f_{kj}^{2}(x)}\leq
$$
$$
\leq C\sum\limits_{k\geq 2}\lambda_{k}^{\frac{d}{2}}\exp(-\frac{\tau}{2}
(\lambda_{k}-\lambda))
\sqrt{\sum\limits_{j=1}^{n_{k}}c_{kj}^{2}}\leq
$$

$$
\leq C\left(\sum\limits_{k\geq 2}\lambda_{k}^{d}
\exp(-\tau (\lambda_{k}-\lambda_{1}))\right)^{\frac{1}{2}}
\left(\sum\limits_{k\geq 2}\sum\limits_{j=1}^{n_{k}}c_{kj}^{2}\right)
^{\frac{1}{2}}
\to 0. \eqno(6)
$$

Combining the condition of the theorem and (6), we conclude that
second summand of right part of (5) converges to zero when
 $\tau\to\infty$.

Now we shall estimate of the right part of inequality (3).

By definition of measure $\nu_{\tau}$ we have
$$
R_{\tau}:=\sum\limits_{k=1}^{N(Q,\tau)}u^{2}(\tau,\gamma_{k})=
\int\limits_{Q}\left[J_{1}(\tau,x)+J_{2}(\tau,x)+J_{3}(\tau,x)\right]
d\nu_{\tau}(x).
$$

Here
$$
J_{1}(\tau,x):=
\exp(-\frac{\tau}{2}\lambda_{1})
\left(\sum\limits_{j=1}^{n_{1}}
c_{1j}f_{1j}(x)\right)^{2};
$$
$$
J_{2}(\tau,x):=2\sum\limits_{j=1}^{n_{1}}c_{1j}f_{1j}(x)
\sum\limits_{k\geq 2}\exp(-\frac{\tau}{2}\lambda_{k})
\sum\limits_{j=1}^{n_{k}}c_{kj}f_{kj}(x);
$$
$$
J_{3}(\tau,x):=
\exp(-\frac{\tau}{2}\lambda_{1})
s_{\tau}^{2}(x).
$$

Observe that we have $F(x)\geq 0, x\in \bar Q$, because
the function $u(t,x)\geq 0$ for all $t\geq 0,~~ x\in  \bar Q$.

Let $M=\max\limits_{x\in \bar Q}F(x)$. Applying (6), we get
 under ~$\tau\to\infty$

$$
|J_{1}(\tau,x)|\leq
\exp(-\frac{\tau}{2}\lambda_{1})
 M^{2}\to 0;
$$

$$
|J_{2}(\tau,x)|\leq 2 M
\exp(-\frac{\tau}{2}\lambda_{1})\sup\limits_{x\in Q}|s_{\tau}(x)|\to 0;
$$

$$
|J_{3}(\tau,x)|\leq
\exp(-\frac{\tau}{2}\lambda_{1})\sup\limits_{x\in Q}s_{\tau}^{2}(x) \to 0.
$$

These inequalities and theorem's condition guarantee the convergence
of ~$R_{\tau}$ to zero under $\tau\to\infty$.

The proof of theorem is complete.

\enddemo

\bigskip

\bigskip

{\bf Example.}
Now we shall investigate particular case of the general problem.
Here we can calculate relevant values of normalizing function and on
the other hand this case may be represent the first approximation
of real situation.

Consider circle domain $Q$ in $E^{2}$:
$x^{2}+y^{2}\leq r^{2}$.
   Assume that the diffusion particles start from point
$(x_{k},y_{k})\in Q$ at the moment $t=0$. The movement of particles
is described by the following stochastic differential equations

$$
d\xi(t)=
\sum\limits_{i}^{2}b_{i}dw_{i}(t)\eqno(7)
$$
$$
\xi(0)=\xi_{0}=(x_{k},y_{k}).
$$

where $b_{1}=(\sigma,0),b_{2}=(0,\sigma)$ and $W(t)=(w_{i}(t),i=1,2)$
be a 2-dimensional Wiener process.

Assume that the equation (7) defines a diffusion process with absorption
on the boundary
 $\partial Q=\{(x,y): x^{2}+y^{2}=r^{2}\}$.

In what follows, the $J_{0}(x), J_{1}(x)$ are Bessel functions zero and first order. It
are defined as the solutions of next equations

$$
\frac{d^{2} y}{d x^{2}} + \frac{1}{x}\frac{d y}{d x} + (1-\frac{n^{2}}{x^{2}})
=0,
$$
$$
y(x_{0})=0,~~~  (x_{0}=\sqrt{\lambda}r);\quad  |y(0)|<\infty;
$$

under $n=0$ and $n=1$.

The value of $\mu_{m}^{(0)}$ is equal to $m$- th
root of equation $J_{0}(\mu)=0$ [6,7].

We will use the  symbol $mes(\cdot)$ to denote the Lebesgue measure.

Put

$$
f(\tau):=\exp\left(-\frac{\tau}{2}
\left(\frac{\sigma \mu_{1}^{(0)}}{r}\right)^{2}
\right).
$$

\bigskip
Suppose $N(\cdot,\tau)$ satisfies the condition
$$
\lim\limits_{\tau\to\infty}N(B,\tau)f(\tau) = mes(B) ,\quad
B\in \Sigma_{mes}(C).
$$

The initial position of $k$-th particle in this case has form
$\gamma_{k}=(x_{k},y_{k})$

The  value of $u(\tau,\gamma_{k})$ is defined as value of $u(t,x,y)$
in this point.
Here $u(t,x,y)$ be a solution of
the following initial-boundary problem

$$
\frac{\partial u}{\partial t}=\frac{\sigma^{2}}{2}\left(\frac{\partial^{2}}
{\partial y^{2}}+\frac{\partial^{2}}{\partial x^{2}}\right)u\quad (x,y)\in Q
,\quad t>0;
$$

$$
u(0,x,y)= 1\quad (x,y)\in Q;
$$

$$
u(t,x,y)=0 \quad \hbox{if}\quad  (x,y)\in \partial D, \quad t\geq 0\eqno(8)
$$

According to results from
book [6, sec. VI] the solution of (8) in the point $(\tau,\gamma_{k})$
can be expressed in form

$$
u(\tau,\gamma_{k})=
\sum\limits_{m=1}^{\infty}
\frac{2 J_{0}\left(\mu_{m}^{(0)}\frac{\sqrt{x_{k}^{2}+y_{k}^{2}}}{r}\right)}
{\mu_{m}^{(0)}J_{1}\left(\mu_{m}^{(0)}\right)}
 \exp\left(-\frac{\tau}{2}
\left(\frac{\sigma\mu_{m}^{(0)}}{r}\right)^{2}\right).
$$

Let us compute the parameter $a$. In this case  it is convenient to
decompose of circle  $C$ by concentric circles
for constuction of integal sums of integral $\int\limits_{C}F(x)\nu(dx)$~~[8].

Define this partition
$$
K_{ni}=\left\{ (x,y)\in C: \frac{ri}{n}\leq  \sqrt{x^{2}+y^{2}}
< \frac{r(i+1)}{n}\right\},\quad 0\leq i\leq n-1.
$$

Now ~$mes( K_{ni})=g(\frac{i+1}{n})-g(\frac{i}{n})$~ where ~
~$g(\rho)=\pi r^{2}\rho^{2}, ~~~0\leq \rho\leq 1$.

Finally, the parameter of Poisson distribution is equal to

$$
a=2\left(\mu_{1}^{(0)}J_{1}(\mu_{1}^{(0)})\right)^{-1}2\pi r^{2}
\int\limits_{0}^{1}J_{0}(\mu_{1}^{(0)}\rho)\rho d\rho=\pi
\left(\frac{2 r}{\mu_{1}^{(0)}}\right)^{2}.
$$

We used the following known relation
~~$\alpha J_{0}(\alpha)=\left[\alpha J_{1}(\alpha)\right]'$ [6, sec.VI]
for calculation of the last integral.

{\bf Acknowledgments.} ~~ The authors are grateful to N.I.Portenko,
A.A. Dorogovtsev and A.M. Kulik for several useful and stimulating discussiones.

\enddocument